\documentclass[12pt,a4paper]{article}
\usepackage{amsmath,supertabular,tabularx,amssymb,amsfonts,float}
\usepackage{amsmath,amsfonts}
\usepackage{graphicx}
\usepackage{subfigure}
\usepackage{epstopdf}
\usepackage{longtable}

\input epsf
\usepackage{epsfig}

\newcommand{\bc}{\begin{center}}
\newcommand{\ec}{\end{center}}
\newcommand{\bq}{\begin{quote}}
\newcommand{\eq}{\end{quote}}
\newcommand{\btab}{\begin{tabular}}
\newcommand{\etab}{\end{tabular}}

\newcommand{\be}{\begin{equation}}
\newcommand{\ee}{\end{equation}}
\newcommand{\beqa}{\begin{eqnarray*}}
\newcommand{\eeqa}{\end{eqnarray*}}
\newcommand{\beqn}{\begin{eqnarray}}
\newcommand{\eeqn}{\end{eqnarray}}
\newcommand{\bbibl}{\begin{thebibliography}{99}}
\newcommand{\ebibl}{\end{thebibliography}}

\newcommand{\nn}{\nonumber}

\newcommand{\ba}{\begin{array}}
\newcommand{\ea}{\end{array}}

\newcommand{\rad}{\stackrel{d}{\rightarrow}}

\newcommand{\qed}{\rule{1.3mm}{3mm}}

\newcounter{cnt1}
\newcounter{cnt2}
\newcounter{cnt3}
\newcommand{\blr}{\begin{list}{$($\roman{cnt1}$)$} {\usecounter{cnt1}
                \setlength{\topsep}{0pt} \setlength{\itemsep}{0pt}}}
\newcommand{\bla}{\begin{list}{$($\alph{cnt2}$)$} {\usecounter{cnt2}
                \setlength{\topsep}{0pt} \setlength{\itemsep}{0pt}}}
\newcommand{\bln}{\begin{list}{$($\arabic{cnt3}$)$} {\usecounter{cnt3}
                \setlength{\topsep}{0pt} \setlength{\itemsep}{0pt}}}
\newcommand{\el}{\end{list}}

\renewcommand{\rad}{\stackrel{d}{\rightarrow}}

\begin{document}
\title{\bf Relative Efficiency of Higher Normed
Estimators Over the Least Squares Estimator} \vspace{1cm}


\author{Gopal K Basak\thanks{ Stat-Math Unit,  Indian Statistical Institute,  Kolkata, India,  email: gkb@isical.ac.in} \;   
Samarjit Das\thanks{Corresponding author:  ERU, Indian Statistical Institute, 203 B.T. Road, Kolkata-700108, India, 
Tel: +91-33-2575-2627, email: samarjit@isical.ac.in} \; 
Arijit De \thanks {Former Masters student of Indian Statistical Institute} \; 
and Atanu Biswas \thanks { Applied Statistics  Unit,  
Indian Statistical Institute,  Kolkata, India,  email: atanu@isical.ac.in} 
\thanks { We are thankful to Probal Choudhury, Jyoti Sarkar for their constructive inputs; 
and also Sattwik Santra for providing us the data.} \\ Indian Statistical Institute \\
}

\date{}
\maketitle

\begin{abstract}
  In this article, we study the performance of the
  estimator that minimizes   $L_{2k}- $ order loss function (for $ k \ge \; 2 )$    against 
  the estimators 
  which minimizes the 
  $L_2-$ order loss function (or the least squares estimator). Commonly occurring examples illustrate the differences 
  in efficiency between  $L_{2k}$  and $L_2 -$ based estimators.  We derive an empirically testable condition under 
  which  the $L_{2k}$  estimator is more efficient than the least squares estimator.  We construct a  simple decision 
  rule to choose between $L_{2k}$ and $L_2$ estimator. Special emphasis is provided to study $L_{4}$  estimator.    
  A detailed simulation study verifies the effectiveness of this decision rule. Also, the superiority of the $L_{2k}$ 
  estimator is demonstrated in a real life data set.  
	
\end{abstract}

\vspace{0.2cm}
\noindent {\bf Key Words:}  Least Squares, Higher Order Loss Function.

\vspace{0.2cm}
\noindent {\bf JEL Classification: C01, C13 }

\newpage

\section{Introduction}
\setcounter{equation}{0}
$\quad$

The least squares (LS hereafter) method is possibly the most popular method of estimation  routinely used to  estimate the 
underlying (regression) parameters.   Stigler (1981) rightly said:  "The method of least squares is the automobile of modern 
statistical analysis: Despite its limitations, occasional accidents, and incidental pollution, it and its numerous variations, 
extensions, and related conveyances carry the bulk of statistical analysis, and are known and valued by nearly all". Such an 
overwhelming popularity of the LS may be due to its simplicity,  optimal properties  and  robustness to any distributional 
assumption. Moreover,  it leads to the best (minimum variance) estimator under normality. Laplace used the name "most 
advantageous method".
However, it appears to us that such an irresistible popularity of the LS may have impeded the exploration of other {\it smooth} 
loss functions. Comparative computational difficulty might be another reason that such exploration  was not favoured by  pioneers 
such as  Gauss, Laplace and others. Whereas, a large literature  to incorporate  {\it non-smooth} loss functions  in order to 
address  (outlier) robustness have been developed. Unfortunately and surprisingly, whole statistical literature is somewhat mute
on the possible use of   {\it smooth} higher order loss functions.  Therefore, it is a pertinent question to ask: Are there any 
relative advantages in using higher order smooth loss function  compared to the omnipresent least squares? 
Our aim here is to study  an appropriate  higher order  estimator  and compare its efficiency against the LS. In the regression 
set up, we find a significantly large and useful class of error distributions for which a higher order loss function is more 
efficient than the LS.  In this paper, we  give an  empirically testable condition under which a higher order smooth loss functions 
lead to a more efficient estimator than the LS. We also provide a simple but pragmatic decision rule  to make a choice between  
$L_{2k}$ and $L_2$. A detailed simulation study shows the effectiveness of such a decision rule.  

In Section 2, we describe the model and develop the methodology needed to compare the efficiency of different loss functions.
 Section 3 provides various classes of error distributions which are used for comparison of estimator. In Section 4, we provide 
 a decision rule along with its asymptotic properties. Section 5 provides an epilogue where we consider  very general classes of
 parametric distributions on finite support to illustrate the enormous scope of applicability of the  $L_4 -$ based loss functions.
Section 6 summarizes the results of simulation study  on mixture distributions. Section 7 gives an application to  real life data.
Section 8 ends with some concluding remarks and identifies possible future directions of research.

\section{Model and assumptions}
\setcounter{equation}{0}

Consider a linear regression set up
\be
\label{gen-regression}
Y = X \beta + \varepsilon ,
\ee
where  $Y$ is an $n \times 1$ vector of observations, $X$ is an $n \times k$ design matrix and $F$ is the cumulative 
distribution function (cdf) corresponding to the error vector the error $\varepsilon. $  We also assume the following 
regular conditions:

\begin{enumerate}

\item $\mbox{plim} \left( \frac{X' X}{n} \right) = A$, a finite and nonsingular matrix.
  
\item $E(\varepsilon | X) = 0$.
  
\item $E(\varepsilon \varepsilon' | X) = \sigma^2 I$.
  
\item Observations are independent.
  
\end{enumerate}

    In this set up, the ordinary least squares estimate (OLS) ${\widehat \beta}_{OLS}$ is the {\it best linear unbiased estimator} 
    in the sense of minimum variance. 
It is well-known that for the OLS estimator, under $\beta =\beta_0$,
$$\sqrt{n} \left( {\widehat \beta}_{OLS} - \beta_0\right) \rad N \left( 0,
\sigma^2 ( X' X)^{- 1}\right) .$$

Furthermore, it is clear that the minimization of $\sum_{i=1}^{n}(Y_i-X_i'\beta )^k$ is pointless if $k$ is odd; and the 
minimization of $\sum_{i=1}^{n}|Y_i-X_i'\beta |$ or $\sum_{i=1}^{n}|Y_i-X_i'\beta |^k$ for odd $k$ will not be very convenient 
because  of lack of differentiability $ wrt \;\; \beta$.  Therefore, it remains to check   whether minimization of 
$\sum_{i=1}^{n}(Y_i-X_i'\beta )^{2k}, $ for some positive integer $k$ other than 1, can yield a better result than 
${\widehat \beta}_{OLS}$, at least in some cases. 
If so, our objective is to identify those cases. It is obvious that the corresponding estimators for such a cases will 
be non-linear. Furthermore,  when the error is normal, then the best linear unbiased estimator is indeed  the
best unbiased estimator. Thus, for normal or near normal
error, ${\widehat \beta}_{OLS}$ always will be better than any other estimator. A closer look reveals that  deviation
from  uni-modality  causes  the robustness properties of LS to falter. 

Studying the efficacy of higher order normed based estimator is important on its own right; not necessarily in comparison to least square. It opens up the possibility to consider a convex combination  of loss functions of various degrees. Arthanari and 
Dodge  (1981) considers convex combination of $L_1$ and $L_2$ norms; and studies its properties. Convex combination
of $L_1, \ldots, L_p$ may lead to more useful estimator; and the resultant estimator is expected to be robust to any
distributional assumption.  Earlier also, such an use of higher order loss functions is attempted.  Turner (1960) heuristically touch upon the possible 
use of a higher order loss function in the context of estimation of the location parameter. He discusses several kinds of 
general PDFs;  and advices in 
the case of the double exponential, to minimize the sum of the absolute
deviations; in the case of the normal, to minimize the sum of the squared
deviations (least squares); and in the case of the q-th power distribution,
to minimize the sum of the q-th  power of the deviations (least q-th's). 
Attempts are also made to define a general class of likelihood to derive a  robust parameter estimates, robust to distributional assumption. For example, Zeckhauser and Thompson (1970)  defines a general class of distribution; and empirically found its suitability.



\subsection{Methodology}

For the exposition purpose, let us first consider  the simple bivariate linear regression model
$$Y = \alpha + \beta x_i + \varepsilon_i,$$
for $i=1,\ldots ,n$. The usual approach to take the error function as $S_2 = \sum^n_{i = 1} ( Y_i - \alpha - \beta x_i)^2$.
We obtain ${\widehat \theta}_{OLS} = ({\widehat \alpha}_{OLS} ~~{\widehat \beta}_{OLS})'$ by minimizing $S_2$ with respect to 
$\alpha$ and $\beta$. Note that ${\widehat \theta}_{OLS}$ is the best estimates in the class of
linear unbiased estimators. Hence there may be some nonlinear estimator with better efficiency.

In contrast, we shall take  $S_4 = \sum^n_{i = 1} ( Y_i - \alpha - \beta x_i)^4$ as our
loss function,  and derive  ${\widehat \theta}_{L_4} = \left( {\widehat \alpha}_{L_4} ~~{\widehat \beta}_{L_4}\right)$
as the estimator of $\theta$. Our objective is to compare ${\widehat \theta}_{OLS}$
and ${\widehat \theta}_{L_4}$,  and discover conditions under which the latter performs better than the former.
Both of these estimators are $M$-estimators. So they possess  the
properties such as consistency and asymptotic normality under some standard conditions.

For the OLS estimator,
\be
\sqrt{n} \left( {\widehat \theta}_{OLS} - \theta_0\right) \rad N \left( 0, \sigma^2 S^{- 1}\right) ,\ee
where
$$ S= \left(\begin{array}{cc}
  n & \sum_{i=1}^{n} x_i\\
  \sum_{i=1}^{n} x_i & \sum_{i=1}^{n}x^2_i
\end{array}\right).$$

We exhibit that ${\widehat \theta}_{L_4}$ satisfies  the following result.

{\bf Lemma 1.} \be \sqrt{n} \left( {\widehat \theta}_{4} - \theta_0\right) \rad N \left( 0, \frac {\mu_6 -\mu_3^2} {9\mu_2^3} S^{- 1}\right).\ee

{\it Proof:}  For the $L_4$ estimator, using the M-estimator property, we have 
$$\sqrt{n} \left( {\widehat \theta}_{L_4} - \theta_0\right) \rad N \left( 0, V (
\theta_0)\right) ,$$
where 
$$V ( \theta_0) = A ( \theta_0)^{- 1} B ( \theta_0) [ A ( \theta_0)^{-
1}]',$$
\beqn
B ( \theta_0) & = & E [ \psi ( y, \theta_0) \psi ( y, \theta_0)'] - E [ \psi ( y, \theta_0)] E[\psi ( y, \theta_0)]',\nn \\
A ( \theta_0) & = & E \left( \frac{\delta}{\delta \theta} \psi (
y, \theta_0) \right),\nn
\eeqn
and 
$$\psi =\frac{\delta S_4} {\delta \theta} = ( - 4) \left(\begin{array}{c}
  \sum ( Y_i - \alpha - \beta x_i)^3\\
  \sum ( Y_i - \alpha - \beta x_i)^3 x_i
\end{array}\right).$$

Let  $\mu_k$ denote the $k$th order central moment corresponding to the distribution of $\varepsilon$.

Consequently, we get
$$
E[ \psi ( y, \theta_0)] = (- 4) \left(\begin{array}{c}
  n \mu_3 \\
  \mu_3 \sum x_i
\end{array}\right) , \ \ \ \mbox{and hence}, \ \ \ E[ \psi ( y, \theta_0)] E[ \psi ( y, \theta_0)]' 
=  16 \mu_3^2 R 
$$
where 
$$
R = \left(\begin{array}{cc}
  n^2  & n  \sum x_i \\
  n \sum x_i & (\sum x_i)^2 \\
\end{array}\right) . 
$$
Also,
$
E [ \psi ( y, \theta_0) \psi ( y, \theta_0)'] = 16 ( \mu_6 S + \mu_3^2 Q) $, 
where
$$Q = \left(\begin{array}{cc}
 n(n-1) & (n-1)\sum_{i=1}^{n}x_i\\
 (n-1)\sum_{i=1}^{n} x_i & (\sum_{i=1}^{n}x_i)^2 - \sum_{i=1}^{n}x^2_i
\end{array}\right) 
$$
Thus,
$B ( \theta_0) = 16 \mu_6 S + 16 \mu_3^2 (Q - R)$.
Since 
\begin{eqnarray*}
Q - R &=& \left(\begin{array}{cc}
 n(n-1) - n^2 & (n-1)\sum_{i=1}^{n}x_i - n \sum_{i=1}^{n}x_i\\
 (n-1)\sum_{i=1}^{n} x_i - n \sum_{i=1}^{n}x_i & (\sum_{i=1}^{n}x_i)^2 - \sum_{i=1}^{n}x^2_i - (\sum_{i=1}^{n}x_i)^2
\end{array}\right) \nn\\
&=&  \left(\begin{array}{cc}
  - n & - \sum_{i=1}^{n}x_i \\
 -  \sum_{i=1}^{n}x_i &  - \sum_{i=1}^{n}x^2_i 
\end{array}\right) \ = \ - S ,
\end{eqnarray*}
$$B ( \theta_0) = 16 S (\mu_6 - \mu_3^2) = 16 S Var(\varepsilon^3).
$$
Similarly, one can simplify  $$A ( \theta_0) = 12 \mu_2 .S$$

 Then
$$V ( \theta_0) = A ( \theta_0)^{- 1} B ( \theta_0) [ A ( \theta_0)^{- 1}]' =
\frac{1}{12 \mu_2} S^{- 1} ( 16 S (\mu_6 - \mu_3^2)) ) \frac{1}{12 \mu_2} S^{-1} 
= \frac{\mu_6 - \mu_3^2}{9 \mu_2^2} S^{- 1}  .$$
Hence the proof. \qed

{\bf Theorem 1.}  The $L_4$ estimator performs better than the OLS estimator in terms of
precision iff
\be
\frac{\mu_6 - \mu_3^2}{ \mu_2^3}   < 9 \label{condition1}
\ee

{\bf Proof.}
Proof follows by comparing 2.2 and 2.3. \qed


For symmetric distribution of $\varepsilon$, or whenever $\mu_3=0$,
 the criterion will be
\begin{equation}
\label{condition}
\frac{\mu_6}{9 \mu_2^3} <1.
\end{equation}

Clearly this condition may or may not be satisfied depending on the distribution of $\varepsilon$.

So far, for exposition purpose, we have dealt with a simple regression framework. Now the scope of this paper demands to present all these above findings in a more general regression set-up. The following remark is made to this end. 

\bigskip
\noindent
{\bf Remark 1.}

For a  multiple linear  regression model  with $k$ regressors,
 all the above calculations can be carried out with  $S = X'X$
where
$$X =
= \left(\begin{array}{ccccc}
  1  & x_{11} & x_{21} & \cdots & x_{k1} \\
  1  & x_{12} & x_{22} & \cdots & x_{k2} \\
  \vdots & \vdots & \vdots & \vdots & \vdots \\
  1  & x_{1n} & x_{2n} & \cdots & x_{kn} \\
\end{array}\right) . 
$$  
Whereas $Q$ matrix would be given by
$$
Q = \left(\begin{array}{ccccc}
 n(n-1)  & (n-1) \sum_{i=1}^n x_{1i} & (n-1) \sum_{i=1}^n x_{2i} &   \cdots & n \sum_{i=1}^n x_{ki} \\
(n-1) \sum_{i=1}^n x_{1i} & \sum_{i\neq j} x_{1i} x_{1j} &  \sum_{i\neq j} x_{1i} x_{2j} & \cdots & \sum_{i\neq j} x_{1i} x_{kj} \\ 
(n-1) \sum_{i=1}^n x_{2i} & \sum_{i\neq j} x_{2i} x_{1j} &  \sum_{i\neq j} x_{2i} x_{2j} & \cdots & \sum_{i\neq j} x_{2i} x_{kj} \\ 
  \vdots & \vdots & \vdots & \vdots & \vdots \\
(n-1) \sum_{i=1}^n x_{ki} & \sum_{i\neq j} x_{ki} x_{1j} &  \sum_{i\neq j} x_{ki} x_{2j} & \cdots & \sum_{i\neq j} x_{ki} x_{kj} \\   
\end{array}\right), 
$$
and the matrix $R$ is given by,
$$
R = 
= \left(\begin{array}{ccccc}
 n^2  & n \sum_{i=1}^n x_{1i} & n \sum_{i=1}^n x_{2i} &   \cdots & n \sum_{i=1}^n x_{ki} \\
n \sum_{i=1}^n x_{1i} & (\sum_{i=1}^n x_{1i})^2 &  (\sum_{i=1}^n x_{1i}) (\sum_{i=1}^n x_{2i}) & \cdots & (\sum_{i\neq j} x_{1i}) (\sum_{i=1}^n x_{ki}) \\ 
n \sum_{i=1}^n x_{2i} & (\sum_{i=1}^n x_{2i}) (\sum_{i=1}^n x_{1i}) x_{1j} &  (\sum_{i=1}^n x_{2i})^2 & \cdots & (\sum_{i=1}^n x_{2i}) (\sum_{i=1}^n x_{2i} x_{ki}) \\ 
  \vdots & \vdots & \vdots & \vdots & \vdots \\
n \sum_{i=1}^n x_{ki} & (\sum_{i=1}^n x_{ki}) (\sum_{i=1}^n x_{1i}) &  (\sum_{i=1}^n x_{ki}) (\sum_{i=1}^n x_{2i}) & \cdots & (\sum_{i=1}^n x_{ki})^2 \\   
\end{array}\right) .
$$
Thus, \ $Q-R = - S$ \ which gives the earlier result
that is, $L_4$ estimators are better than that of $L_2$
iff  
$$ 
\frac{\mu_6 - \mu_3^2}{9 \mu_2^3}   < 1 .
$$

\bigskip

It may be interesting to examine the  performance of  $L_{2k}$ relative to that of  $L_2$, or  that of $L_{2k-2}.$  The following two corollaries are presented to this end. 

{\bf Corollary 1.}  $L_{2k}$ estimator better than LS (i.e., $L_2$)  iff
$$\frac{ Var(\varepsilon^{2k-1})}{(2k-1)^2 (Var(\varepsilon))^{2k-1}} < 1 .
$$

{\bf Proof.} The proof is analogous to that of Theorem 1. \qed

{\bf Corollary 2.} $L_{2k}$ estimator better than  $L_{2k-2}$  iff

$$\frac{ Var(\varepsilon^{2k-1})  (2k-3)^2}{(2k-1)^2 (Var(\varepsilon^{2k-3}) ) \mu_2 ^2} < 1 . $$

{\bf Proof.} The proof is analogous to that of Theorem 1. \qed






\section{OLS versus $L_4$ for some selected  distributions}

In this section, we consider few important parametric error distributions  to illustrate the vast scope of applicability of $L_4$ based loss function. The list of distributions considered is no way exhaustive, but certainly shows the immense opportunity of applications in diverse areas.  Now, we check the aforementioned condition (2.4) hold for different distributions of $\varepsilon$.

\subsubsection{U-Shaped  Distribution }
Consider a simple U-shaped distribution: $$f(x)=d x^{2k}; \;\;\; -c \le x \le c, \; \mbox{where}  \; k \mbox{ is a positive integer}.$$
Note that $d=\frac{2k+1}{2c^{2k+1}}.$
It is easy to calculate $$\frac{\mu_6}{9 \mu_2^3}=\frac{(2k+3)^3} {9(2k+1)^2(2k+7)}.$$
Note that for $k=1, \; \frac{\mu_6}{ 9 \mu_2^3}=\frac{125}{729}, $ and in the limit,  $ \frac{\mu_6}{ 9 \mu_2^3}=\frac{1}{9}.$


 A U-shaped distribution has two modes; and can be looked upon as a mixture of two (J-shaped) distributions - a  mixture of  an extreme positively skewed and  another extreme negatively skewed distributions. 
 One popular applied example  of a U-shaped distribution is   the number of deaths  at various ages. Several more  examples can be found in by B. S. Everitt ( 2005).

\subsubsection{Uniform($-a,a$)}

This is  a symmetric distribution. Here
$$\mu_r = \frac{a^{r + 1} - ( - a)^{r + 1}}{( r + 1) \{ a - ( - a)\} },$$
and consequently
$$\frac{\mu_6}{9 \mu_2^3} = \frac{3}{7} < 1.$$
Hence, $L_4$ estimator is better than the OLS estimator, when the error component has  uniform distribution.

\subsubsection{Normal($\mu, \sigma^2$)}

This is again a symmetric distribution where
$$\mu_{2 r} = \sigma^{2 r} ( 2 r - 1) \times ( 2 r - 3) \times \cdots
\times 5 \times 3 \times 1,$$
and hence
$$\frac{\mu_6}{9 \mu_2^3} = \frac{15}{9} > 1.$$
Hence, for  normally  distributed errors, the OLS estimator is  always preferred over the $L_4$ estimator.

\subsubsection{Laplace($\lambda$)}

Here we have
$$\mu_r = \lambda^r \Gamma ( r + 1),$$
if $r$ is even, and, consequently,
$$\frac{\mu_6}{9 \mu_2^3} = \frac{6!}{72} > 1.$$
Hence, when $\varepsilon$ follows Laplace distribution,  the OLS estimator is preferred over the $L_4$ estimator.

\subsubsection{Beta($a,b$)}

The beta distribution is a family of continuous probability distributions defined on the interval [0, 1] parametrized by two 
positive shape parameters, denoted by $a$ and $b$, that appear as exponents of the random variable and control the shape of 
the distribution.  This class of distributions include a  variety of symmetric, 
bell-shaped, positively skewed, negatively skewed, uniform, and 'U-shaped' distributions.   The general form of the central moments 
of the beta distribution are quite
complicated. So we will start with the raw moments and obtain the forms of
$\mu_6$, $\mu_3$ and $\mu_2$. Note that, here
$$\mu_r' = \prod^{r - 1}_{i = 0} \left\{ \frac{a + i}{a + b + i}\right\} .$$
The figure 4 depicted in Section 5, provides a huge range of parameters  for which $L_4$ is  better than $L_2.$




\subsubsection{Gaussian mixture distribution}

Suppose  that $F ( x) = \frac{1}{2} N (\xi_1,
\sigma^2) + \frac{1}{2} N (\xi_2, \sigma^2)$. We assume, for simplicity,  common $\sigma^2$  for both the components.  Here
$$\mu_r = \frac{1}{2} \sum_{i = 0}^r \binom{r}{i} ( \xi_1 - \xi )^{r -
i} \mu_{1i} + \frac{1}{2} \sum_{i = 0}^r \binom{r}{i} ( \xi_2 -
\xi )^{r - i} \mu_{2i},$$
where $\xi = \frac{1}{2} ( \xi_1 + \xi_2)$ and $\mu_{ji}$ is the $i$th central moment of $N(\xi_j,\sigma^2)$ distribution, $j=1,2$. Let $\xi_2 - \xi = \frac{1}{2} ( \xi_2 - \xi_1) = - (\xi_1 - \xi ) = a_0$. Then (\ref{condition}) reduces to
\begin{equation}
\label{mixture-condition}
6 + 18 c^2 - 12 c^4 - 8 c^6 < 0, 
\end{equation}
where $c = \frac{\xi_2 - \xi_1}{2 \sigma}$. The left hand side expression of (\ref{mixture-condition}), which is a function of $c$,
is plotted in Figure 1. From the plot we can see for $|c|>1.058$, this function assumes values less
than zero and then it decreases rapidly. This means that the condition (\ref{condition}) of superiority of $L_4$ estimators will be satisfied if the
means of the two component of the mixture distribution are more than $1.058\sigma$ distance. A mixture of
more than two Gaussian distributions will behave similarly with respect to this condition. Also the case of unequal $\sigma^2$ can be tackled similarly.

\begin{figure}[H]\label{mixture-gauss}
 \centering
\includegraphics[scale=0.4]{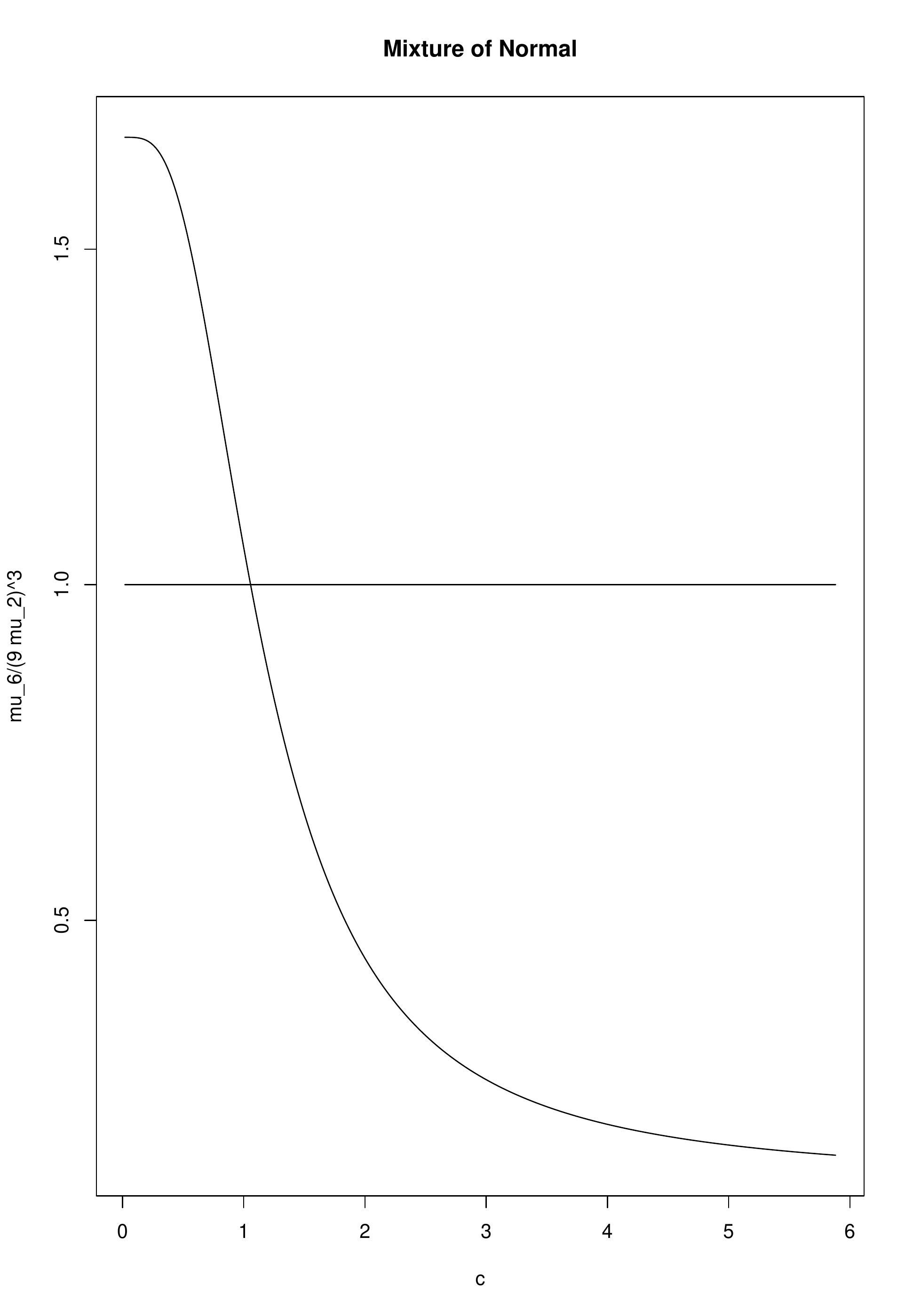}
 \caption{line plot as a function of c}
\end{figure}

\subsubsection{Truncated Normal distribution}
Consider (for simplicity) the both side truncated standard normal distribution. The  even order moments are

$$ \mu_{2k} ^{c}=\frac{2}{d} \int_{0} ^{c}  \frac {  x^{2k} exp(- \frac{x^2}{2}) }  {\sqrt{2 \pi}} dx.  $$
Define $\Delta= \sqrt{2 \pi} (\Phi(c)-0.5).$
Therefore, $$ \mu_{2k} ^{c}=\frac {  c^{2k-1} exp(- \frac{c^2}{2}) }{\Delta}   +\mu_{2k-2} ^{c}  (2k-1) .$$
Now it is easy to  calculate $ \mu_{6} ^{c} =15- \frac{c^5  exp(- \frac{c^2}{2})}{\Delta} -\frac{ 5 c^3  exp(- \frac{c^2}{2})}{\Delta} - \frac{ 15c \; exp(- \frac{c^2}{2})}{\Delta},$ and\\
$ \mu_{4} ^{c} =3- \frac{c^3  exp(- \frac{c^2}{2})}{\Delta} -\frac{ 3c \; exp(- \frac{c^2}{2})}{\Delta},$ and \\
$ \mu_{2} ^{c} =1- \frac{c\;   exp(- \frac{c^2}{2})}{\Delta}. $
Now one can see that as long as $c\in(0, 2.33), \; L_4$ performs better than $L_2.$
One implication of this result is that 97 percent times $L_4$ performs better than $L_2. $

\begin{figure}[H]\label{Truncated Normal}
 \centering
 \includegraphics[scale=0.4]{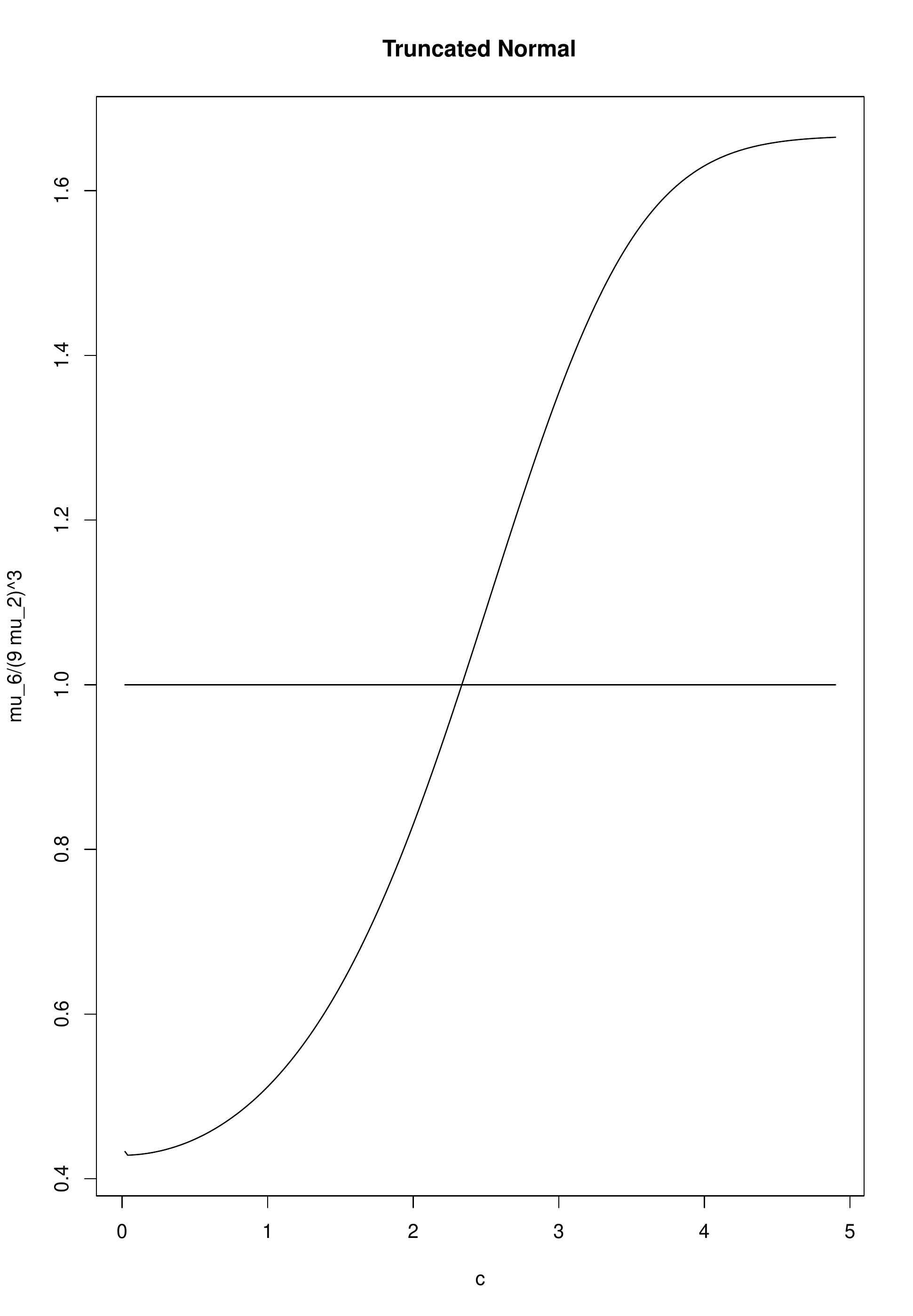}
 \caption{Truncated Normal plot}
\end{figure}

\subsubsection{Raised cosine distribution}

If $X$ follows a raised cosine distribution with parameters $a$ and $b$, denoted by $X\sim COR(a,b)$, then the probability density 
function (pdf) is given by
$$f(x) = \frac{1}{2b}\left[ 1+\cos\left(\pi\frac{x-a}{b}\right)\right], ~~a-b\leq x\leq a+b, ~a\in\mathbb{R}, ~b>0.$$
The form of this distribution resembles that of a normal distribution except for the fact that it has finite tails. Suppose it can be assumed that the value of systematic errors lies in some known interval; and manufacturer has aimed to make device as accurate as possible. In such circumstances, Raised Cosine distribution may be appropriate. Another  popular application is in circular data. See Rinne (2010, pp. 116). Other properties like the cdf, moment generating function (mgf), characteristic functions, 
raw moments up to order 4, and the kurtosis are available in Rinne (2010, pp. 116-118). It is observed that the distribution 
has a kurtosis of 2.1938, less than that of normal distribution. It has a thin tail. Here, using the mgf, we have
\beqn
\mu_6 & = & \frac{b^6 \left(\pi^6 -42\pi^4 + 840\pi^2 - 5040\right)}{7\pi^6},\nn \\
\mu^2 & = & \frac{b^2 \left(\pi^2 - 6\right)}{3\pi^2},\nn
\eeqn
and hence
$$\mu_6/(9\mu_2^3) = \frac{3}{7} - \frac{72\pi^4 - 2196\pi^2 +14472}{7(\pi^2 - 6)^3} = 0.8926<1.$$
Hence, for this distribution, (\ref{condition}) is satisfied, and consequently $L_4$ is preferred for parameter estimation.

\subsection{A Sub-Gaussian family of distributions}

Sub-Gaussian family of distributions is a well-studied family of distribution whose tail is dominated by the normal distribution. As we observed that the $L_4$ estimators are preferred for a distribution for which the tail is thinner than that of a normal
distribution, here we discuss about a relatively uncommon distribution and validity
of the condition with respect to this distribution. Consider a distribution with pdf of the form
$$f ( x) = c \exp (- x^{2 k}),$$
where $k$ is an integer and $c$ is the normalizing constant, which gives
$$c =  \frac{k}{\Gamma \left( \frac{1}{2 k} \right)}.$$
For  various values of $k$, the pdf of the distribution is drawn in Figure 3. Note that $k=1$ provides the normal curve. As $k$ becomes larger and larger, tail of the distribution tend to collapse. For extremely large $k$, the distribution resembles a symmetric curve in a finite support.   It is interesting to consider the peaks of all drawn curves. The first plot (the density plot) of this panel shows that $L_4$ performs better than $L_2$ for all those curves for which peaks are below the red curve. Here it may be mentioned that the red curve is drawn for $k=1.45.$ For the second plot of the panel, various values of $k$  are given in the $x-$axis; and values of the test statistic are given in the $y-$axis. The parallel line, parallel to $x-$axis, shows the cut-off point, which is 1 . The second plot of the panel shows that when the value of $k$ is greater than 1.45, $L_4$ performs better than $L_2$. 

\begin{figure}[H]\label{pdf-k1}
 \centering
 \includegraphics[scale=0.8]{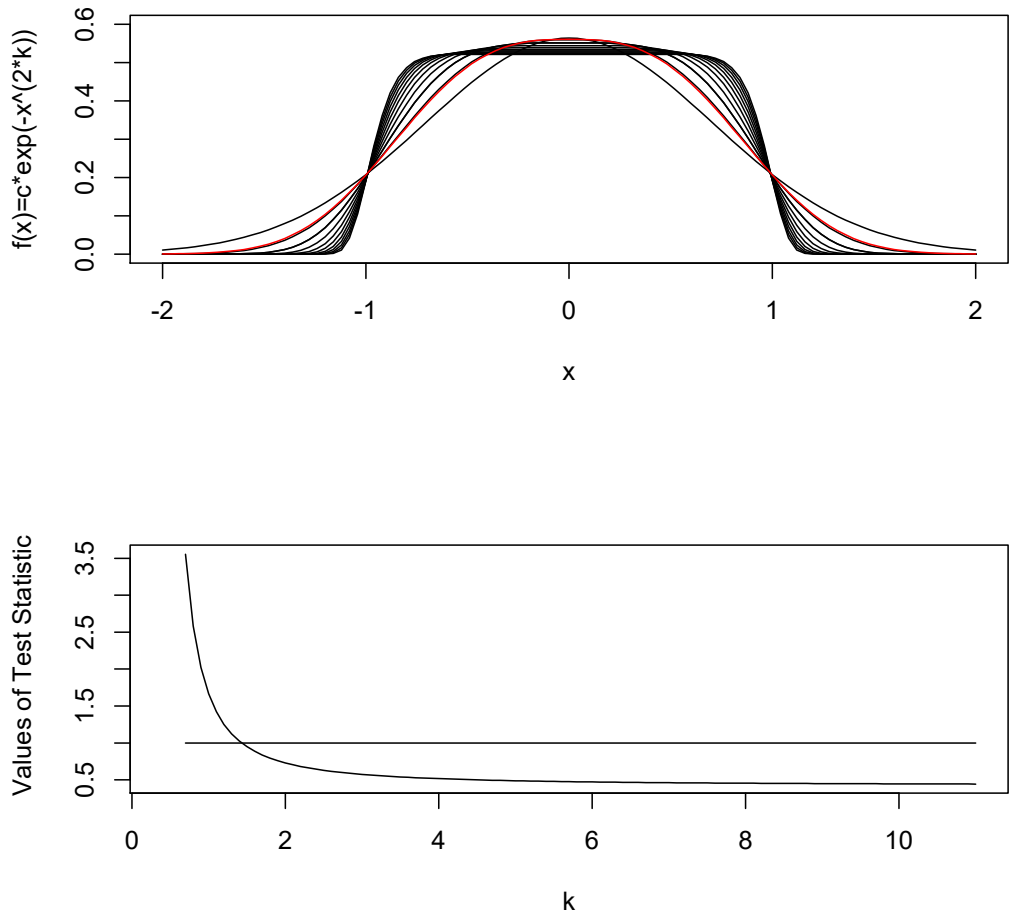}
  \caption{The shape of the various Sub-Gaussian distributions 
  	($ f ( x) = c \exp (- x^{2 k})$ for a range of values $k$ where $L_4$ is better than $L_2$. The red curve is for $k=1.45$. }

\end{figure}

Here $\mu_r = 0$ if r is odd. When r is even, we have
$$\mu_r = c \int_{- \infty}^{\infty} x^r \exp (- x^{2 k}) =  \frac{c}{k}  \Gamma \left(
\frac{r + 1}{k} \right) .$$
We immediately get
$$\frac{\mu_6}{9 \mu_2^3} = \frac{\left( \Gamma \left( \frac{1}{2 k} \right)
\right)^2 \Gamma \left( \frac{7}{2 k} \right)}{9 \left( \Gamma \left(
\frac{3}{2 k} \right) \right)^3}.$$
The values of the test statistic against  various values of $K$ is drawn in the bottom part in Figure 3. We observe that for $k \geqslant 2$, 
$\frac{\mu_6}{9 \mu_2^3}$ assumes values less than 1.

Now, according to  Zeckhauser and Thompson (1970); Turner (1960) and Box and Tiao (1962 ), for a distribution 
with pdf
$$f ( u) = k ( \sigma, m) \exp \left(- \left| \frac{u}{\sigma}
\right|^{m}\right), ~~\sigma > 0, ~m > 0,$$
where $k ( \sigma, m) = \left[ 2 \sigma \Gamma \left(
1 + \frac{1}{m} \right) \right]^{- 1}$, $L_{m}$ estimators dominate all $L_{m'}$ estimators where $m'
< m$.

Note that, for $m = 4$ and $\sigma = 1$, we obtain the distribution displayed in Figure 3. So $L_4$ estimator 
performs better than the corresponding $L_2$ estimator, that is the OLS estimator, which is entirely in agreement
to what we derived above. For the normal
distribution $m = 2 $; $m = 1$ gives the double exponential
distribution; where $m$ tends to $\infty$, the
distribution tends to the rectangular.  The article of Zeckhauser and Thompson (1970) considers  four empirical 
examples  to find that there is a sizable gains in likelihood if $m$ is estimated rather than pre-specified 
equal to 2. All of the evidence they found  leads them to the conclusion that if accurate estimation of a linear
regression line is important, it will
usually be desirable to estimate not only the
coefficients of the regression line, but also the
parameters of the power distribution that generated
the errors about the regression line. The
effect on the estimates of regression coefficients
may not be small.

In the next two section we will construct a decision rule based on the condition (\ref{condition1}) and carry
out some simulation study.

\section{Decision rule: OLS versus $L_4$}

In this Section we derive a decision rule based on the criterion from Section
2 to decide whether OLS or $L_4$ estimator is preferred for some data.\\

{\bf Lemma 2.} Suppose $X$ follows a distribution for which $\mu_r$ exists for all $r$. Then
$$\sqrt{n} \widehat{\mu_r} = \frac{1}{\sqrt{n}} \sum^n_{i = 1} ( x_i -
\overline{x})^r = \frac{1}{\sqrt{n}} \sum^n_{i = 1} ( x_i - \mu)^r - r
\mu_{r - 1} \frac{1}{\sqrt{n}} \sum^n_{i = 1} ( x_i - \mu) + o_p ( 1).$$

{\it Proof:} Observe that
\beqn
 &    & \sqrt{n} \widehat{\mu_r} = \frac{1}{\sqrt{n}} \sum_{i = 1}^n ( x_i -
  \overline{x})^r = \frac{1}{\sqrt{n}} \sum_{i = 1}^n ( x_i - \mu
  + \mu - \overline{x})^r\nn \\
 & = & \frac{1}{\sqrt{n}} \sum_{i = 1}^n ( x_i - \mu)^r + r \frac{1}{n} \sum_{i = 1}^n ( x_i - \mu)^{r - 1} \sqrt{n} ( \mu - \overline{x})\nn \\
 &    & + \binom{r}{2} \frac{1}{n} \sum_{i = 1}^n (  x_i - \mu)^{r - 2} \sqrt{n} ( \mu - \overline{x})^2 + \cdots .\label{expan}
\eeqn
Now all the terms other than the first and second term of (\ref{expan}) are of the order $o_p ( 1)$ because $\sqrt{n} ( \mu - \overline{x})$ is $O_p ( 1)$, and hence $\sqrt{n} ( \mu - \overline{x})^k, ~k \geq 2$, is $o_p( 1)$. 
Also $\frac{1}{n} \sum_{i = 1}^n ( x_i - \mu)^{r - 1} = \mu_{r - 1} + o_p ( 1)$. Therefore
$$\sqrt{n} \widehat{\mu_r} = \frac{1}{\sqrt{n}} \sum^n_{i = 1} (
  x_i - \overline{x})^r = \frac{1}{\sqrt{n}} \sum^n_{i = 1} ( x_i - \mu)^r - r
  \mu_{r - 1} \frac{1}{\sqrt{n}} \sum^n_{i = 1} ( x_i - \mu) + o_p ( 1).$$ \hfill{$\Box$}

Furthermore, by delta method,
\[ \sqrt{n} \left( ( {\widehat \sigma}^2)^{\frac{r}{2}} - ( \sigma^2)^{\frac{r}{2}}
   \right) = \frac{r}{2} ( \sigma^2)^{\frac{r}{2} - 1} \sqrt{n} (
   {\widehat \sigma}^2 - \sigma^2) + o_p ( 1) .\]
Then, we have the following Theorem.\\

{\bf Theorem 2.} Let $v = \frac{\mu_6 - \mu_3^2}{^{} \sigma^6}$. Suppose $\mu_{12}$ exists
  for distribution of $X$. Then
  \[ \sqrt{n} ( {\widehat v} - v) = \frac{\alpha_0}{{\widehat \sigma}^6}
     \frac{1}{\sqrt{n}} \sum_{i = 1}^n Z_i + o_p ( 1) ,\]
  where $\alpha_0 = \left( 1, ~~- (6  \mu_5 - 3 \mu_2 \mu_3^{}), ~~-\mu_3 ,~~- 3 \sigma^4 v\right)$ and 
  $Z_i = \left(\begin{array}{c}
    ( x_i - \mu)^6 - \mu_6\\
    ( x_i - \mu)\\
    ( x_i - \mu)^3 - \mu_3\\
    ( x_i - \mu)^2 - \sigma^2
  \end{array}\right) .$
  
\bigskip  
Proof of this theorem is given in the Appendix B.

\section{An Epilogue}
In this section we consider a very general class of parametric distributions on finite support (this assumption 
is made to ease plot drawing) to illustrate the enormous scope of applicability of $L_4$ based loss function. 
Consider the class of distribution: 
$$ f=d(1+x^2)^a, \ a\in R \ d > 0, \ |x| \le 1, \ . $$ Here $d$ depends on $a$ to make the $f$ a density. 
The first plot of the panel depicts the shape of the density for different values of 
$a$. Depending on the value of $a$, this class of distributions includes various 'U-shaped' (for $a > 0$) 
and 'Bell-Shaped' 
(for $a < 0$) distributions.  It is interesting to consider the peaks of all drawn curves. The first plot (the density plot) of this panel shows that $L_4$ performs better than $L_2$ for all those curves for which peaks are below the red curve. Here it may be mentioned that the red curve is drawn for $a=-3.2.$ For the second plot of the panel, various values of $a$  are given in the $x-$axis; and values of the test statistic are given in the $y-$axis.  The second plot of the panel shows that when the value of $a$ is greater than -3.2, $L_4$ performs better than $L_2$, in all shape going from low deep to hump till it reaches certain level.

\begin{figure}[H]\label{f-plot-surface1}
  \centering
 \includegraphics[scale=0.8]{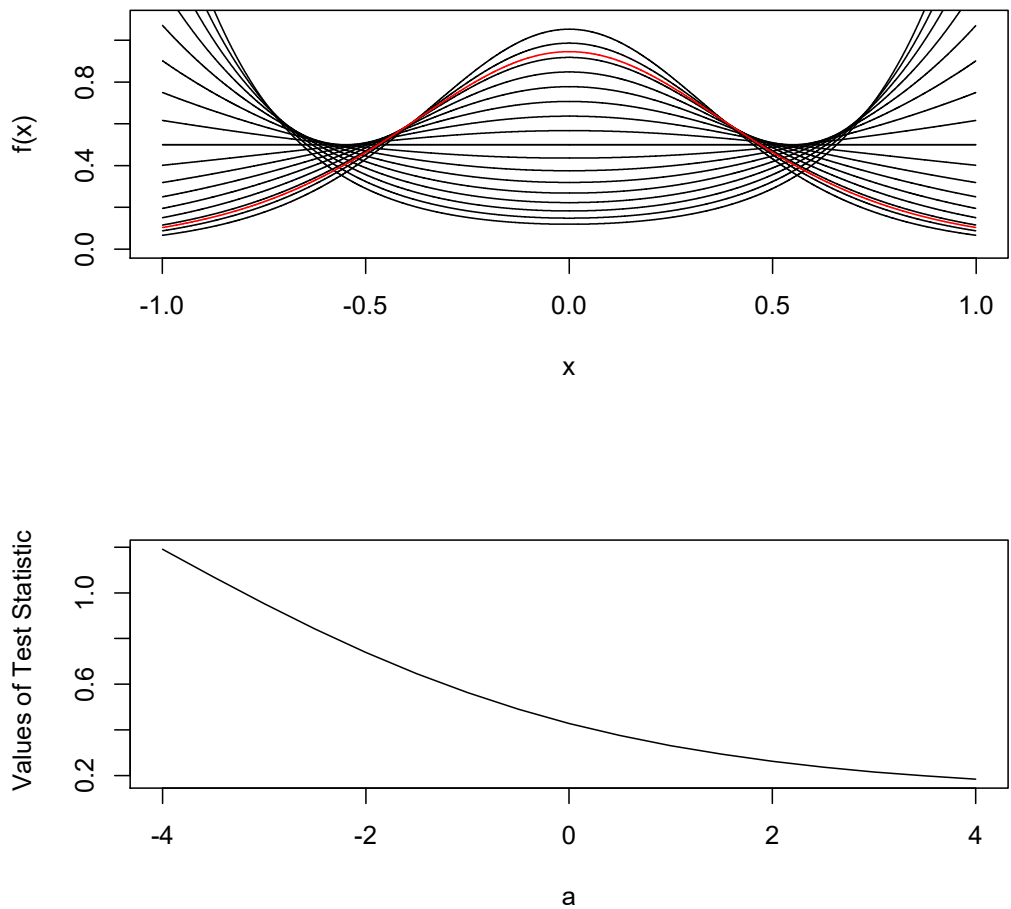}
 \caption{The shape of the distributions 
  ($f=d(1+x^2)^a, \ a\in R, \ d > 0, \ |x| \le 1$) for a range of values $a$ where $L_4$ is better than $L_2$. }
\end{figure}

Plots of another parametric family of distributions (belongs to the 
Pearsonian family, Type II), given by
$$f=d(1-x^2)^a, \ a > - 1 \ d > 0, \ \mbox{ function of } a, \ |x| < 1$$ are shown below. 
  It may be noted that this particular distribution is linked to $Beta$ distribution as well. To see this, let $Y \sim Beta(\alpha_1, \alpha_2), Y \in (0,1).$ Let $X= a+(b-a)Y.$ Then  
$$ f(x)=\frac {\Gamma(\alpha_1+ \alpha_2)} {\Gamma(\alpha_1)  \Gamma(\alpha_2)} \frac {(x-a)^{\alpha_1 -1} (b-x)^{\alpha_2 -1}} { (b-a)^{\alpha_1 +\alpha_2 -1}}.$$ To see the equivalence,  set $ a=-b=-1, $ and  $ \alpha_1=\alpha_2= \alpha +1.$

\begin{figure}[H]\label{f-plot-surface2}
 \centering
 \includegraphics[scale=0.8]{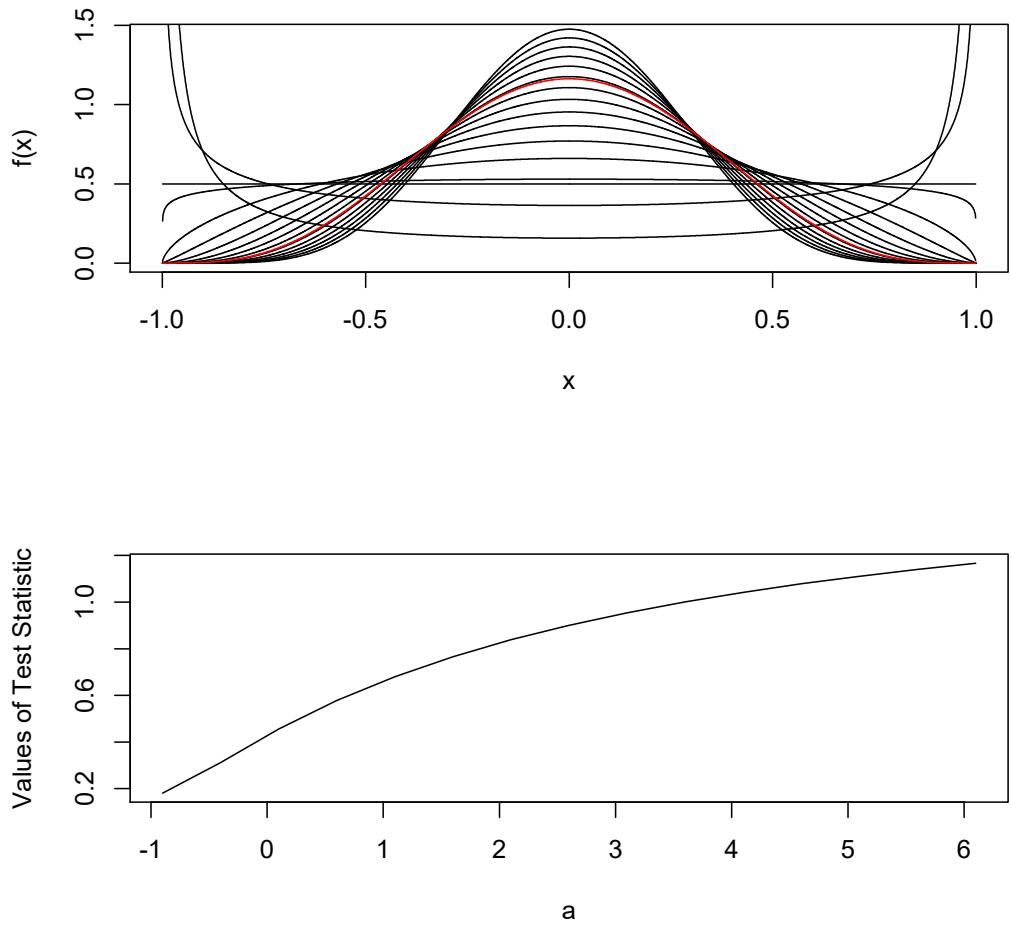}
 \caption{The shape of the distributions ($f=d(1-x^2)^a, \ a > - 1,  \ d > 0, \ |x| \le 1$) for a range of values $a$ where $L_4$ is better than $L_2$. }
\end{figure}

Figure 5 also depicts the same feature as in Figure 4.  Depending on the value of $a$, this group of parametric of distributions 
depicts various 'U-shaped'  and 'Bell-Shaped'. It is interesting to consider the peaks of all drawn curves. The first plot (the density plot) of this panel shows that $L_4$ performs better than $L_2$ for all those curves for which peaks are below the red curve. The red curve is drawn for $a=3.5.$   The second plot of the panel shows that when the value of $a$ is greater than 3.5, $L_4$ performs better than $L_2$, in all shape going from low deep to hump till it reaches certain level.

These two distributions illustrate the enormous possibility of use of  $L_4$ based loss function. Future
studies will investigate whether this is a general
phenomenon for other Pearsonian family of errors distributions. 

\section{Simulation study}
\setcounter{equation}{0}
$\quad$
We carry out the decision making procedure under 0-1 loss function and calculate the risk function, which is the expected loss. 
Here we generate data from three types of distribution, one for which $L_4$ is always better than OLS estimator, one where OLS 
is better than $L_4$, and the third one is near the boundary. The values of the calculated risk are given in Table \ref{table 1}.

\bigskip
\bc
{\em Table 1 given in Appendix A should be here}
\ec
\bigskip
 
Simulation study is based on  10000 iterations; and  with sample sizes of 100,  200, 500,  1000, 2000, 5000.   

The first panel of Table 1 is based on mixture of  two $T$ distributions with 6 degrees of freedom (DF) each. Mean of each 
components are set at $(5, -5), (4,- 4), (3, -3),  (2, -2). $ Here it may be mentioned that our test needs existence of 6th 
order moments. To this end, we need  $t$ distribution with at least 7 df. DF 6 is considered to examine the performance of our 
decision rule even when moments do not exist. Mixture coefficients are taken from $U(0,1)$ distribution.  From this part of the 
table, it is clear that decision is more certain as sample size increases; more importantly, it is so when the distance between 
the two components are more. 

The second  panel of the Table is based on  mixture of two $T$ distribution with 10 df each. Here findings corroborate with the 
first panel. The third panel is based on mixture of two $T$ distribution with 20 df each.

The 4th panel of the table is based on mixture of  two  asymmetric $Beta$ distributions. Mixture coefficients are taken from 
$U(0,1)$ distribution. The fifth panel is based on two symmetric beta distributions with weight from $U(0,1)$. The first column, 
in Panel 5 needs special attention. The parameter combination $ ((4,4; 4,4))$ is chosen such that it is in the neighborhood of 
the boundary the test statistic. Here it shows that test does not favour (for the large sample case, n=5000) any one, as expected.
Here the risk is near 50\% . 

The 6th and 7th panel of the Table are based on mixture of two normal distributions. Here also findings are on the expected line. 
As sample size increases, test correctly discriminates between $L_4$ and $L_2.$


\section{Empirical Illustration}
In this sub-section, we provide two illustrations. One is based on a constructed data set which resembles many real life scenario; and the second one is based on a real life data set. 

\subsection{ Constructed Example}
Data often contains rounding errors. Variables (like heights or weights, age in years, or birth weight in ounces.) that by
their very nature are continuous are, nevertheless, typically measured in a discrete manner. People feel more comfortable to report their age as mid forty, mid fifty and so on. They are rounded to a certain level of accuracy, often to some preassigned
decimal point of a measuring scale (e.g., to multiples of 10 cm, 1 cm, or 0.1 cm) or
simply our preference of some numbers over other numbers.  The reason may be the
avoidance of costs associated with a fine measurement or the imprecise nature of
the measuring instrument. The German military, for example, measures the height
of recruits to the nearest 1 cm. Even if precise measurements are available, they
are sometimes recorded in a coarsened way in order to preserve confidentiality or to
compress the data into an easy to grasp frequency table. 

Here we consider the linear regression where the dependent variable is rounded to nearest integer; independent 
variables are free of any such errors. The dependent variable is generated as
$$ y_{st}= 8 + 1 \times x_1+2 \times x_2  \ \mbox{ where} \ x_1=1.3 \times sample.int(10); x_2=2.32*sample(10:18).$$ 

However, assume that we do not observe $Y_{st}$ but observe $$Y= 5\times floor(y_{st}/5) 
\footnote{sample.int(10) randomly arranged 1 to 10 integers; sample(10:18) randomly arranged 10 to 18 integers; 
floor ($y_{st}/5$) is the largest integer  less than or equal to  $y_{st}$. }   .$$

Now we are regressing $Y$ on $X_1$ and $X_2.$  For this example, we set a moderate sample size of 40. We 
consider 5000 replication.The output is summarized as follows: 


\bc
{\bf Table 2:} Average Estimates  Based on the Constructed Data.\\
\begin{tabular}{c|ccc}
 
\hline
   & \multicolumn{3}{c}{ Parameters} \\ \cline{2-4}
Estimates & Intercept=5.5 & Slope 1 =1 & Slope 2=2\\ \hline
Average ($L_2$)  & 7.013 & 1.115 & 1.924 \\
Average  ($L_4$) & 6.548 & 1.021 & 1.962 \\
 \hline
\end{tabular}
\ec
It is observed that 90 percent times $L_4$ is preferred over $L_2$ based on our proposed decision rule.

After estimation of the model, it may be of interest to know which set of estimators provides the best fit. In the present context it is a tricky problem to find an appropriate 'goodness of fit' measure. Likelihood based methods are not tenable. Similarly, residual sum of square or $R^2$ are not not useful to compare the performances of these two set of parameter estimates. Here we suggest to apply the idea of Pseudo $R^2$  (See Cameron and Trivedi, 2005 for details, page No. 311).  

 Let $Q(\theta)$ denotes the objective function being maximized, $Q_0$
denotes its value in the intercept-only model, $Q_{fit}$ denotes the value in the fitted model,
and $Q_{max}$ denotes the largest possible value of $ Q(\theta)$. Then the maximum potential
gain in the objective function resulting from inclusion of regressors is $Q_{max} - Q_0$ and
the actual gain is $Q_{fit} - Q_0$. This suggests the measure
$$R_{RG}^2= \frac{Q_{fit} - Q_0} {Q_{max} - Q_0}.$$

where the subscript RG means relative gain. Note that, for least squares, $R^2=R_{RG}^2.$ For both the loss functions, $Q_{max}=0.$ 

We also calculated the number of times the Pseudo $R^2$  for $L_4$ is numerically greater than that of $L_2$. It is astonishing  
to see that 100 percent times  the Pseudo $R^2$  for $L_4$ is numerically greater than that of $L_2$.

\subsection{Real Life Example}  

For our empirical analysis, we use the data provided by the National Sample Survey Organization of India viz. the NSSO 68th 
round all India unit level survey on consumption expenditure (Schedule1.0, Type 1 and 2) conducted during July 2011 to June 2012. 
This dataset is a nationally representative sample of household and individual characteristics based on a stratified sampling 
of households. For this round, the dataset is comprised of 1,68,880 household level observations. The dataset provides a 
detailed list of various household and individual specific characteristics along with the consumption expenditures of the 
households. In addition to this, data is also provided on the households' localization which includes the sector (Rural or Urban), 
the district and the state/union territory (henceforth, the union territories will be referred to as states). For our analysis, 
we use the amount of land possessed (in logarithm form) by the households as our principal (dependent) variable along side 
various demographic variables as controls (independent variables). The kernel density plot clearly suggest that amount of land 
possession by rural households  does have a bimodal distribution \footnote{The same phenomenon is also seen for all-India 
households (rural and urban together). The plot presented here is for rural household excluding the households with no land. 
It is interesting to note that bi-modality is observed both for (1) households with non-zero amount of land ; and (2) with all 
households. All the results presented here are based on rural household with non-zero lands. Number of rural households with 
non-zero amount of land is 98483. Whole study is based on This set of 98483 observations.}. The plot clearly indicates that India is suffering from  "vanishing middle-class 
syndrome," only marginal  and rich  farmers are there. The empirical analysis demands some routine and rudimentary  summary 
statistics as given in Table 2.   We regress the amount of land possessed ($Y$)   on six explanatory variables, \footnote{ We also tried with many other explanatory variables available in our master file.  We also repeated the same exercise for all-india (rural and urban together) households. It is need less to mention that overall findings are same across all models we attempted.} viz, Median 
age of a household (mage), the number of children below 15 years of age (chlt15), the number of old people above 60 years of 
age (Ogt60),  the number of male member in the households (male), the number of female member in the households (female); and 
finally the number of member with education level above 10th standard (highedu).   We then estimate \footnote{In this paper we 
do not pursue the endogeneity issue, if any.} the linear regression model based on $L_2$ and $L_4.$ The estimated results are 
summarized as below: 

\begin{figure}[H]\label{Kernel Density}
	\centering
	\includegraphics[scale=0.3]{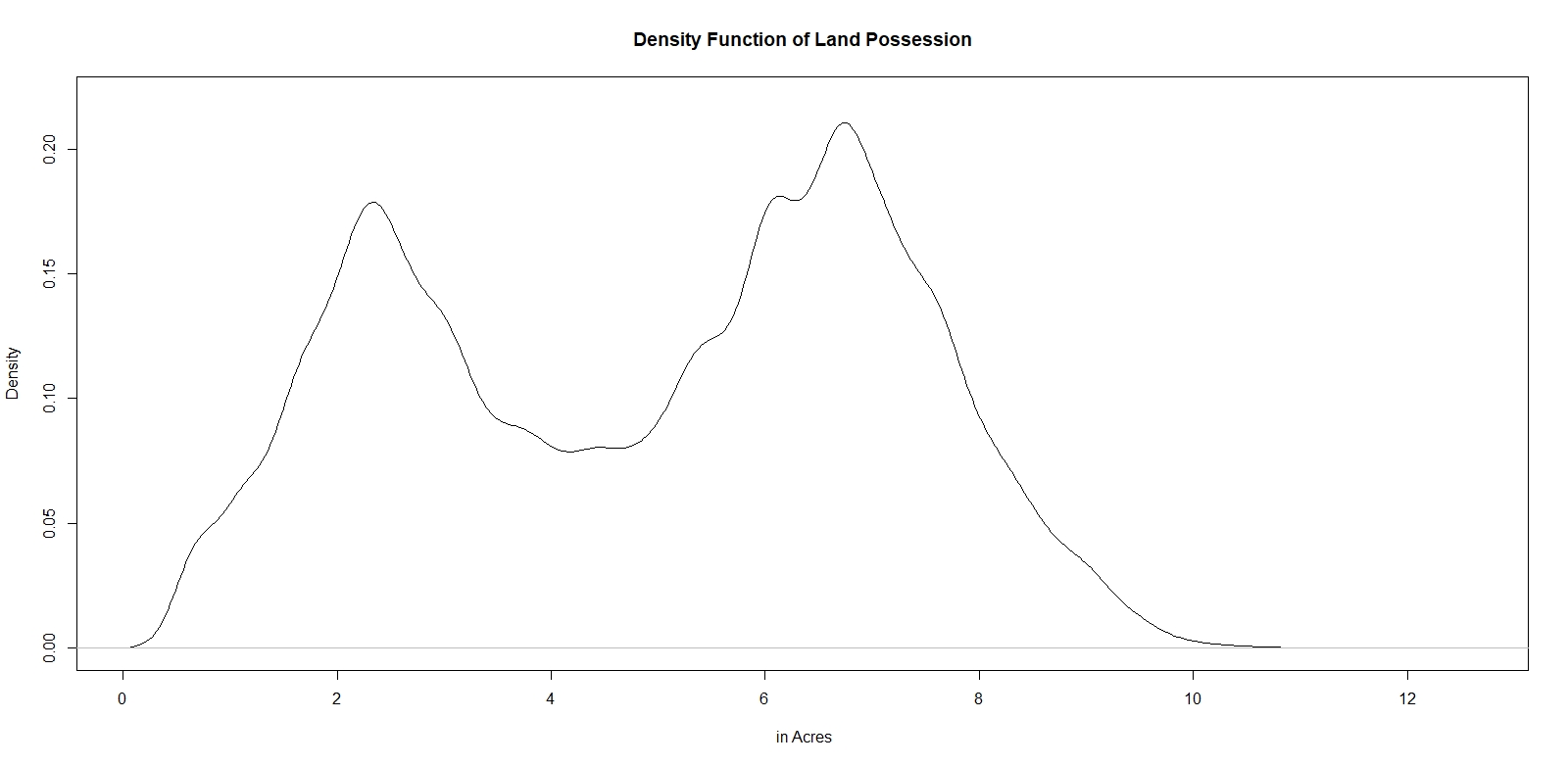}
	\caption{Kernel density plot of amount of land possession by rural households}
\end{figure}


\bc
{\bf Table 3:} Summary Statistics.\\
\begin{tabular}{cccccc}
	
	\hline\hline
	Mean  & Median &  SE & Min & Max & Kurtosis  \\
	\hline
	\hline 
	4.945 & 5.352 &2.290 &0.693 & 12.007 & 1.944 \\ \hline
	\multicolumn{6}{p{15.5cm}}{\footnotesize {\bf Note:} (i) All the results presented here are based on rural household with non-zero lands. Number of rural households with 
		non-zero amount of land is 98483. (ii) This table is based on non-logarithm data. }	
	
\end{tabular}
\ec

\newpage

\bc
{\bf Table 4:} Model Estimates and Standard Errors .\\
	\begin{tabular}{c c c  } 
		\hline\hline 
	Variables &	$L_2$ & $L_4$ \\ 	\hline 
	Intercept &	3.03493710 & 3.47252382  \\ 
		& (0.0342302942 ) & (0.0120790850 )\\
	mage &	0.01240624 & 0.01077902  \\
	& 	(0.0007863559) & ( 0.0002774869) \\
   chlt15 &    -0.14759963  & -0.09232877  \\
     &  (0.0082999771 )  & (0.0029288714 ) \\
	Ogt60 & 	0.06754135 & 0.04168707 \\
	&	( 0.0135886351) &  (0.0047951174 ) \\
	male &	0.36720544 & 0.24403590 \\
	&	(0.0064274676 )  &  (0.0022681058 ) \\
	female &	0.31067372 & 0.21197564 \\
	&	(0.0071582759  ) &  (0.0025259912 ) \\
	highedu &   0.18669815 &  0.16207843 \\
	  &  (0.0074214316 )  & (0.0026188528 ) \\
		\hline 
	\multicolumn{3}{p{15.5cm}}{\footnotesize {\bf Note:} (i) All the results presented here are based on rural household with non-zero lands. Number of rural households with 
		non-zero amount of land is 98483. Whole study is based on this set of 98483 observations. Logarithm transformation is taken 
		to reduce degree of heteroscedasticity. (ii) Standard errors are provided in the parenthesis. (iii) Least-squares based estimates are used as an initial estimates for $L_4$ estimation. (iv) 'Rootsolver' in R is used to obtain the $L_4$ estimates. }	
		
	\end{tabular}
	\ec

	It can be noted that Standard errors (SE) of the parameter estimates are provided in the parenthesis. 
	It is to observe, as expected, that SE of $L_4$ based estimates are significantly and uniformly  less 
	than that of $L_2$ based estimates. The value of our proposed test statistic is 5.30311871  which lies 
	beyond 95 per cent confidence interval (8.90603838 9.09396162) suggesting that  $L_4$ based estimates 
	are more efficient than that of $L_2$ . The pseudo $R^2$ for  $L_4$  is 0.14779372
	 and the same for $L_2$ is 0.09794736 .   The pseudo $R^2$  also clearly suggests the supremacy  of 
	 $L_4$ over $L_2. $
	
\section{Discussion}
\setcounter{equation}{0}
$\quad$
This paper tried to give answer to the unassailable question: Does higher order loss function based estimator 
perform better than the omnipresent least squares? Every teacher, student faces this question on the first-day  
class on regression analysis. We tried to show that, in several real life situations, {\it smooth} higher order 
loss function based estimator may lead to more efficient estimator as compared to universal least squares. It is 
true that least squares has one unassailable advantages, its simplicity. It may also be computationally less intensive. However, with the advent of modern computing power, computational issues may hardly be relevant.

Further work may commence in the following directions.   A generalized version of the condition similar to the 
one derived in section 4 may be useful for comparing $L_{2 k}$ and $L_{2 k + 2}$ estimators. This may be  obtained 
following a similar approach i.e. by obtaining the variance of these two estimators using the expression of variance 
of m-estimators and comparing them.  However, estimation of higher moments may have impact on the performance of the 
proposed decision rule.  It may be useful to study the impact of outliers on the parameter estimates coming from 
higher order based loss function. Comparison of {\em break-down point}  of $L_2,L_4$ estimators may be very useful. It 
may also be interesting to find robust standard errors for $L_4$  for non $i.i.d$ set up.

It may be extremely useful to consider a convex combination  of loss functions of various degrees. Arthanari and 
Dodge  (1981) considers convex combination of $L_1$ and $L_2$ norms; and studies its properties. Convex combination
of $L_1, \ldots, L_p$ may lead to more useful estimator; and the resultant estimator is expected to be robust to any
distributional assumption. Such combination may give   answer to the omnipresent question: What is the optimal loss 
function for a given data set? The choice and design of loss functions is important in any practical application (see 
Hennig and Kutlukaya, 2007).   Future research will shed light in this direction.

\newpage
{\bf Appendix A}
\begin{center}
\begin{longtable}{|l|l|l|l|l|l|}
	\caption{Simulated Risk}
	\label{table 1}\\
	
		\hline
		Sample size & Mixture of T-distribution             & (5, -5)      & (4,- 4)    & (3, -3)    & (2, -2)    \\ \hline
		\endfirsthead
		
		\hline
		\multicolumn{6}{|l|}{Continued to next page...}\\ \hline
		\endfoot
		
		\hline
		\endlastfoot
		
		\hline
		\multicolumn{6}{|l|}{...Continued from previous page} \\ \hline
		\endhead
		
		100         & Mixture of T-distn (df=6)             & 9940         & 9834       & 9200       & 5896       \\
		200         &                                       & 9952         & 9838       & 9256       & 5571       \\
		500         &                                       & 9961         & 9844       & 9309       & 4861       \\
		1000        &                                       & 9965         & 9856       & 9337       & 4309       \\
		2000        &                                       & 9975         & 9896       & 9348       & 3763       \\
		5000        &                                       & 9973         & 9912       & 9445       & 2833       \\ \hline
		100         & Mixture of  T-distn (df=10)           & 10000        & 9986       & 9860       & 8092       \\
		200         &                                       & 9999         & 9992       & 9900       & 8329       \\
		500         &                                       & 10000        & 9997       & 9959       & 8842       \\
		1000        &                                       & 10000        & 9996       & 9973       & 9041       \\
		2000        &                                       & 10000        & 10000      & 9991       & 9368       \\
		5000        &                                       & 10000        & 9999       & 9988       & 9666       \\ \hline
		100         & Mixture of T-distn (df=20)            & 10000        & 9999       & 9981       & 9155       \\
		200         &                                       & 10000        & 10000      & 9996       & 9585       \\
		500         &                                       & 10000        & 10000      & 10000      & 9846       \\
		1000        &                                       & 10000        & 10000      & 10000      & 9950       \\
		2000        &                                       & 10000        & 10000      & 10000      & 9984       \\
		5000        &                                       & 10000        & 10000      & 10000      & 9999       \\ \hline
		& Mixture of  Beta-distn (asym)         & (4,10; 10,4) & (1,4; 4,1) & (2,4; 4,2) & (3,4; 4,3) \\ \hline
		100         &                                       & 10000        & 10000      & 9691       & 3632       \\
		200         &                                       & 10000        & 10000      & 9992       & 4791       \\
		500         &                                       & 10000        & 10000      & 10000      & 7277       \\
		1000        &                                       & 10000        & 10000      & 10000      & 9241       \\
		2000        &                                       & 10000        & 10000      & 10000      & 9957       \\
		5000        &                                       & 10000        & 10000      & 10000      & 10000      \\ \hline
		& Mixture of   Beta-distn (sym)         & (4,4; 4,4)   & (3,3; 3,3) & (2,2; 2,2) & (1,1; 1,1) \\ \hline
		100         &                                       & 2007         & 3740       & 7665       & 9987       \\
		200         &                                       & 1881         & 4748       & 9445       & 10000      \\
		500         &                                       & 2089         & 7178       & 9993       & 10000      \\
		1000        &                                       & 2574         & 9200       & 10000      & 10000      \\
		2000        &                                       & 3464         & 9957       & 10000      & 10000      \\
		5000        &                                       & 5810         & 10000      & 10000      & 10000      \\ \hline
		& Mixture of two normal distributions   & (3, -3)      & (2, -2)    & (1, -1)    & (0, 0)     \\ \hline
		100         &                                       & 9999         & 9806       & 1862       & 211        \\
		200         &                                       & 10000        & 9959       & 1330       & 26         \\
		500         &                                       & 10000        & 10000      & 815        & 0          \\
		1000        &                                       & 10000        & 10000      & 478        & 0          \\
		2000        &                                       & 10000        & 10000      & 265        & 0          \\
		5000        &                                       & 10000        & 10000      & 63         & 0          \\ \hline
		& Mixture of three normal distributions & (4, -4)      & (3, -3)    & (2, -2)    & (1, -1)    \\ \hline
		100         &                                       & 9968         & 9560       & 5589       & 550        \\
		200         &                                       & 10000        & 9951       & 6589       & 181        \\
		500         &                                       & 10000        & 10000      & 8252       & 17         \\
		1000        &                                       & 10000        & 10000      & 9450       & 1          \\
		2000        &                                       & 10000        & 10000      & 9922       & 0          \\
		5000        &                                       & 10000        & 10000      & 10000      & 0        \\ \hline 

\end{longtable}
\end{center}

\bigskip

{\bf Appendix B}
  
{\it Proof of Theorem 2:} We write
\beqn
{\widehat v} - v & = & \frac{\widehat{\mu_6} - \widehat{\mu_3}^2}{{\widehat \sigma}^6} - \frac{\mu_6 -
  \mu_3^2}{\sigma^6}= \frac{{\widehat \mu}_6 - {\widehat \mu}_3^2 - \mu_6 +
  \mu_3^2}{{\widehat \sigma}^6} - v \frac{{\widehat \sigma}^6 -
  \sigma^6}{{\widehat \sigma}^6}\nn \\
 & = & \left[ \frac{\frac{1}{n} \sum_{i = 1}^n ( x_i - \overline{x})^6 - \left(
  \frac{1}{n} \sum_{i = 1}^n ( x_i - \overline{x})^3 \right)^2 - \mu_6 +
  \mu_3^2}{{\widehat \sigma}^6} \right] - v \left[ \frac{(
  {\widehat \sigma}^2)^3 - ( \sigma^2)^3}{{\widehat \sigma}^6}
  \right],\nn
\eeqn
and hence
\beqn
 &  & \sqrt{n} ( {\widehat v} - v)\nn \\
 & = & \frac{1}{\hat{\sigma}^6} \left[ \frac{1}{\sqrt{n}} \sum_{i = 1}^n ( x_i -
  \overline{x})^6 - \left( \frac{1}{\sqrt{n}} \sum_{i = 1}^n ( x_i -
  \overline{x})^3 \right) \left( \frac{1}{n} \sum_{i = 1}^n ( x_i -
  \overline{x})^3 \right) - \sqrt{n} ( \mu_6 - \mu_3^2) \right)\nn \\
 &    & - v \sqrt{n} \left[ \frac{( {\widehat \sigma}^2)^3 - (
  \sigma^2)^3}{{\widehat \sigma}^6} \right]\nn \\
 & = & \frac{1}{{\widehat \sigma}^6} \left[ \frac{1}{\sqrt{n}} \sum_{i = 1}^n (
  x_i - \mu)^6 - 6 \mu_5 \frac{1}{\sqrt{n}} \sum_{i = 1}^n ( x_i - \mu) -
  \left( \frac{1}{\sqrt{n}} \sum^n_{i = 1} ( x_i - \mu)^3 - 3 \mu_2
  \frac{1}{\sqrt{n}} \sum^n_{i = 1} ( x_i - \mu) \right) \mu_3 \right.\nn \\
 &    & \left. - \sqrt{n} ( \mu_6 - \mu_3^2) - 3 \sigma^4 v \sqrt{n} (
  {\widehat \sigma}^2 - \sigma^2) \right] + o_p ( 1)\nn \\
 & = & \frac{1}{{\widehat \sigma}^6} \left[ \frac{1}{\sqrt{n}} \sum_{i = 1}^n ( ( x_i
  - \mu)^6 - \mu_6) - (6 \mu_5 - 3 \mu_2 \mu_3) \frac{1}{\sqrt{n}} \sum_{i =
  1}^n ( x_i - \mu) - \mu_3 \left( \frac{1}{\sqrt{n}} \sum^n_{i = 1} ( ( x_i -
  \mu)^3 - \mu_3) \right) \right.\nn \\
 &    & - \left. 3 \sigma^4 v \left( \frac{1}{\sqrt{n}} \sum_{i = 1}^n
  ( ( x_i - \mu)^2 - \sigma^2) \right) \right] + o_{p ( 1)}\nn \\
 & = & \frac{1}{\hat{\sigma}^6} \left(\begin{array}{cccc}
    1, & - ( 6 \mu_5 - 3 \mu_2 \mu_3), & -\mu_3 & - 3 \sigma^4 v
  \end{array}\right) \frac{1}{\sqrt{n}} \left(\begin{array}{c}
    \sum_{} [ ( x_i - \mu)^6 - \mu_6]\\
    \sum ( x_i - \mu)\\
    \sum_{} [ ( x_i - \mu)^3 - \mu_3]\\
    \sum_{} [ ( x_i - \mu)^2 - \sigma^2]
  \end{array}\right) + o_p ( 1)\nn \\
 & = & \frac{\alpha_0}{{\widehat \sigma}^6} \frac{1}{\sqrt{n}} \sum_{i = 1}^n Z_i + o_p ( 1).\nn
\eeqn
\hfill{$\Box$}

Now,
$$\frac{1}{\sqrt{n}} \sum_{i = 1}^n Z_i = \sqrt{n} \overline{Z}
\rad N ( 0, \Gamma),$$
where $\Gamma = \lim_{n \rightarrow
\infty} nE ( \overline{Z}  \overline{Z}')$ and $\overline{Z} =
\frac{1}{n} \sum_{i = 1}^n Z_i$. Here
$$\overline{Z}= \left(\begin{array}{c}
  m_6^o - \mu_6\\
  m_1^o\\
  m_3^o - \mu_3\\
  m_2^o - \sigma^2
\end{array}\right),$$
with
$$m_r^o = \frac{1}{n} \sum_{i = 1}^n ( x_i - \mu)^r$$
and
\beqn
E ( m_r^o) & = & \mu_r, \nn \\
V ( m_r^o) & = & \frac{\mu_{2 r} - \mu_r^2}{n}, \nn \\
Cov( m_r^o, m_s^o) & = & \frac{\mu_{r + s} - \mu_r \mu_s}{n},\nn
\eeqn
and hence
$$\Gamma = \lim_{n \rightarrow \infty} 
{nE} ( \overline{Z}  \overline{Z}') =
\left(\begin{array}{cccc}
  \mu_{12} - \mu_6^2  & \mu_7 & \mu_9 - \mu_3 \mu_6 & \mu_8 - \mu_2 \mu_6\\
  \mu_7 & \mu_2 & \mu_4 & \mu_3\\
  \mu_9 - \mu_3 \mu_6 & \mu_4 & \mu_6 - \mu_3^2 & \mu_5 - \mu_2 \mu_3\\
  \mu_8 - \mu_2 \mu_6 & \mu_3 & \mu_5 - \mu_2 \mu_3 & \mu_4 - \mu_2^2
\end{array}\right).$$

Consequently, we get
$$\sqrt{n} ( {\widehat v} - v) \rad N \left( 0,
\frac{\alpha_0 \Gamma \alpha_0'}{{\sigma}^{12}} \right),$$
where
$$s^2 = \widehat{V ( v)} = \frac{{\widehat \alpha}_0 {\widehat \Gamma}
{\widehat \alpha}_0'}{{\widehat \sigma}^{12}}.$$

Hence our decision theoretic problem reduces to
$$A_0 : v \geq 9 ~~~~\mbox{ versus } ~~~~A_1 : v < 9.$$
We use the statistic
$$T= \frac{\sqrt{n} ( {\widehat v} - 9)}{s},$$
which is asymptotically $N(0,1)$, under $A_0$. \qed


\section{Reference}

\begin{enumerate}

\item Arthanari, T. S. and Dodge, Y. (1981), “Mathematical Programming in Statistics,” New York:
John Wiley


\item Everitt B. S. (2005) U-Shaped Distribution, Encyclopedia of Biostatistics, 8, John Wiley \& Sons.  



\item Box, G. E. P., and G. C. Tiao (1962): A Further Look at Robustness via Bayes's Theorem, {\it  Biometrika },
49, 419-432.

\item Cameron A. C. and P. Trivedi (2005): Microeconometrics, Cambridge.  

\item Hennig. C, and M. Kutlukaya( 2007): Some  Thoughts about the Design  of
Loss Functions, {\it REVSTAT – Statistical Journal}
 5,  1,  19–39.


\item Rinne, H. (2010). Location-Scale Distributions -- Linear Estimation and Probability Plotting Using MATLAB. {\it Preprint.}

\item  Stigler. S., (1981): Gauss and the Invention  of Least  Square, {\it The Annals of Statistics}, 9,3, 465-474.

\item Turner, M. C. (1960): On Heuristic Estimation Methods, {\it Biometrics} , 16,  299-301.

\item  Zeckhauser, R.,  and M.  Thompson (1970): Linear Regression with Non-Normal Error Terms. {\it The Review of Economics and Statistics}, Vol. 52, No. 3, 280-286

\end{enumerate}

\end{document}